\documentclass{article}
\usepackage[dvips]{graphicx}
\usepackage{amsthm}
\usepackage{amsfonts}
\usepackage[dvips]{hyperref}

\newenvironment{definition}{\paragraph{Definition.}}{}
\newtheorem{lemma}{Lemma}
\newtheorem{theorem}[lemma]{Theorem}
\newtheorem{corollary}[lemma]{Corollary}

\newcommand{\Sinfty}{S_\infty}
\newcommand{\ch}{\mathrm{ch}}
\newcommand{\nfrac} [2] { \displaystyle \frac{#1}{#2} }
\renewcommand{\mid}{|}
\renewcommand{\triangle}{\mbox{\large$\bigtriangleup$}}
\newcommand{\seg}[2]{\ensuremath{#1\hspace{-0.1em}#2}}

\providecommand{\R}{\ensuremath{\mathbb{R}}}
\newcommand{\N}{\ensuremath{\mathbb{N}}}
\newcommand{\Q}{\ensuremath{\mathbb{Q}}}

\newcommand{\A}{\ensuremath{{\cal A}}}

\newcommand{\Saltinfty}{\tilde{S}_\infty}
\newcommand{\Salt}{\tilde{S}}
\newcommand{\triangleboundary}{{\cal T}}
\renewcommand{\line}{\mathrm{line}}

\title{On the Density of Iterated Line Segment Intersections}%

\author{
  Ansgar Gr{\"u}ne\footnotemark[1] \hspace{1cm}
  Sanaz Kamali Sarvestani\footnotemark[1]\\[1cm]
  Department of Computer Science I,
  University of Bonn,
  Germany
}
\date{December 2005, revised December 2006 and May 2007}

\begin{document}

\maketitle
\renewcommand{\thefootnote}{\fnsymbol{footnote}}
\footnotetext[1]{E-Mail: {\tt \{gruene,kamali\}@cs.uni-bonn.de}}
\renewcommand{\thefootnote}{\arabic{footnote}}

\begin{abstract}\noindent
  Given~$S_1$, a finite set of points in the plane, we define a
  sequence of point sets~${S_i}$ as follows:
  With~$S_i$ already
  determined, let~$L_i$ be the set of all the line segments connecting
  pairs of points of~$\bigcup_{j=1}^i S_j$, and let $S_{i+1}$ be the
  set of intersection points of those line segments in $L_i$, which
  cross but do not overlap.
  We show that with the exception of some starting configurations
  the set of all crossing points~$\bigcup_{i=1}^\infty S_i$
  is dense in a particular
  subset of the plane with nonempty interior.
  This region is the intersection of all closed half planes which contain
  all but at most one point from~$S_1$.

\end{abstract}


\section{Introduction}\label{sec:intro}

Given $S = S_1$, a finite set of points in the Euclidean plane,
let~$L_1$ denote the set of line segments connecting pairs of points
from $S_1$.
Next, let $S_2$ be the set of all the intersection
points of those line segments in~$L_1$ which do not overlap.
We continue to define sets of line segments~$L_i$ and point sets~$S_i$
inductively by
\begin{eqnarray*}
  L_i
  & := &
  \bigg\{
    pq
    \; \bigg| \;
    p,q \in \bigcup_{j=1}^{i}S_j
    \;\wedge\;\;
    p\neq q
  \bigg\}
  ,
  \\ \nonumber
  S_{i+1}
  & := &
  \left\{
    x
    \;\mid\;
    \{x\} = l\cap l' \mbox{ where }
    l, l' \in L_i
  \right\}
  .
\end{eqnarray*}
Finally, let $S_{\infty}:=\bigcup_{i=1}^\infty S_i$ denote the limit
set.
\begin{figure}[htbp]
  \begin{center}
    \includegraphics{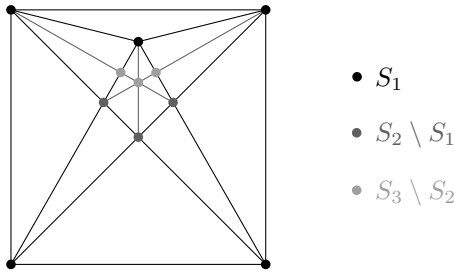}
    \caption{\label{fig:Example} The first iterations $S_1$, $S_2$ and $S_3$
      of segment intersections.}
  \end{center}
\end{figure}

In this article we show which starting configurations~$S$ give rise to
density of~$S_\infty$.
And for those, we prove that $S_{\infty}$
is dense in a certain region~$K(S)$.


Several results concerning the density of similar iterated constructions
are known.
Bezdek and Pach~\cite{Pach} studied the following problem:
Let $S_1 = \{O_1, O_2\}$ consist of two distinct points in the plane with
distance less than~$2$.
Define $S_{i+1}$ inductively as the set of all intersection points
between the unit circles, whose centers are in $S_i$.
They proved,
if the distance between~$O_1$ and~$O_2$ does not equal~$1$ nor~$\sqrt{3}$,
the limit point set~$\bigcup_{i=1}^\infty S_i$
is everywhere dense in the plane.

Kazdan~\cite{Kaz} considered two rigid motions $A$ and $B$
applied to an arbitrary point~$p$ in the plane.
An implication of his more powerful result is that
the set $\Sinfty := \bigcup_{i=1}^{\infty} S_i$
defined by the iterations~$S_1 := \{p\}$
and $S_{i+1}$ $:=$ $\{Ax,$ $A^{-1}x,$ $Bx,$ $B^{-1}x \mid x \in S_i\}$
is everywhere dense in the plane
if $A$ is a rotation through an angle incommensurable with~$2\pi$
and $A$ and $B$ do not commute.

B{\'a}r{\'a}ny, Frankl, and Maehara~\cite{Frankl} showed that with the exception
of four special cases, the set of vertices of those triangles
which are obtained from a particular starting triangle~$T$
by repeated edge-reflection is everywhere dense in the plane.


Ismailescu and Radoi\v{c}i\'{c}~\cite{Ismailescu} examined a
question very similar to ours.
The only difference is that they considered lines instead of line segments.
They proved by applying nice elementary methods
that with the exception of two cases the crossing points
are dense in the whole plane.
Hillar and Rhea~\cite{Hillar} independently proved the same statement with
different methods.

Despite the similarity of the two settings,
the methods from the analysis of iterated line intersections
cannot be transferred to our case of segment intersections.
The latter problem turns out to be more difficult.
It has more exceptional cases, where the crossing points
are not dense in any set with non-empty interior.
And in non-exceptional cases the crossing points
are not dense in the whole plane but only
in a particular convex region~$K(S)$.
In fact, if we do not consider exceptional configurations,
the line result is a direct consequence
of the segment result presented here.
If the iterated intersections of line segments
are dense in a region with non-empty interior,
clearly,
the intersections of lines starting with the same point set
are dense in the same region and thereby in the whole plane.

Our work is motivated by another interesting problem introduced
recently by Ebbers-Baumann et al.~\cite{Klein};
namely how to embed a given finite point set into
a geometric graph of small dilation.
Here a {\em geometric graph} is a graph in the
Euclidean plane, where the vertices are points in the plane,
the edges are rectifiable curves connecting the two adjacent vertices,
and the edge lengths equal the lengths of the corresponding curves.
Given such a geometric graph~$G$,
for any two vertices~$p$ and~$q$
we define their \emph{vertex-to-vertex dilation} as
\begin{displaymath}
  \delta _G(p,q) := \nfrac{|\pi(p,q)|}{|pq|},
\end{displaymath}
where $\pi(p,q)$ is a shortest path from~$p$ to~$q$ in~$G$ and
$|\,.\,|$ denotes the Euclidean length.
The \emph{dilation of $G$} is defined by
\begin{displaymath}
  \delta(G)
  :=
  \sup_{\mbox{\scriptsize $p$, $q$  vertices of  $G$},\;p\neq q}
  \delta _G(p,q).
\end{displaymath}
A geometric graph~$G$ of smallest possible dilation $\delta(G)=1$
is called \emph{dilation-free}.
We will give a list of all cases of dilation-free planar graphs
in Section~\ref{sec:definitions}.
It can also be found at Eppstein's geometry junkyard~\cite{Junk}.
Given a point set~$S$ in the plane, the \emph{dilation of $S$} is defined by
\begin{displaymath}
  \Delta (S)
  :=
  \inf
  \left\{
    \;\delta(G)\;
    \left|\;
      G=(V,E) \mbox{ planar graph},\;
      S \subseteq V
      \mbox{ and } V \setminus S \mbox{ finite }
    \right\}
  \right.
  .
\end{displaymath}
Determining $\Delta(S)$ seems to be very difficult.
The answer is even unknown if $S$ is a set of five points
placed evenly on a circle.
However, Ebbers-Baumann et al.~\cite{Klein} were able to prove
$\Delta(S)\leq 1.1247$ for every finite point set $S \subset \R^2$
and they showed lower bounds for some special cases.

A natural idea for embedding $S$ in a planar graph of small dilation
is to try to find a geometric graph~$G=(V,E)$, $S\subseteq V$,
such that $\delta_G(p,q)=1$ for every $p,q\in V$.
Now, suppose we have found such $G$.
Obviously, for every pair $p,q\in S$, the line segment $pq$ must be
a part of $G$.
Since $G$ must be planar, every intersection point of these line segments
must also be in $V$ and so on.

If this iteration produces only finitely many intersection points,
i.e.\ the set~$\Sinfty$ defined at the very beginning of this article
is finite,
we have a planar graph~$G=(V,E)$ with $S \subseteq V$,
$|V\setminus S|<\infty$ and
$\delta(G)=1$; thus $\Delta(S)= 1$.
This shows that $|\Sinfty|<\infty$ can only hold if $S$ is
a subset of the vertices of a dilation-free planar graph.
We call those point sets {\em exceptional configurations}.
Note that $\Delta(S)=1$ could still hold for other sets.
There could be a sequence of proper geometric graphs
whose dilation does not equal $1$ but converges to $1$.

\begin{figure}[htb]%
  \begin{center}%
    \includegraphics{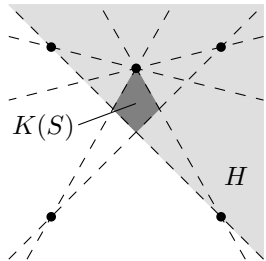}%
    \caption{\label{fig:Candidate} The candidate~$K(S)$
      and a closed half plane~$H$
      containing all but one starting point.}%
  \end{center}%
\end{figure}

Our main result is the following.
\begin{theorem}\label{satz:theo}
  Let $S=S_1$ be a set of $n$ points in the plane, which is not an
  exceptional configuration.
  Then $\Sinfty$ is dense in the region~$K(S)$.
\end{theorem}
The region~$K(S)$ can be described by the following simple
definition.
An example is shown in Figure~\ref{fig:Candidate}.
Until we prove that $\Sinfty$ is indeed dense in~$K(S)$,
we call~$K(S)$ the candidate.
\begin{definition}\label{def:Candidate}
The \emph{candidate} $K(S)$ is defined as the
intersection of those closed half planes which contain all the
starting points except for at most one point,
\begin{displaymath}
K(S):= \bigcap_{p \in S}\,\bigcap_{H\supset S\setminus\{p\}}H.
\end{displaymath}
The second intersection is taken over all closed half planes $H$
which contain $S\setminus\{p\}$.
\end{definition}

After the first version of our proof was worked out,
Klein and Kutz~\cite{Kutz}
used the statement of Theorem~\ref{satz:theo}
to prove the first non-trivial lower bound $\Delta(S) \geq 1.0000047$
which holds for every non-exceptional finite point set~$S$.
They also sketched a different proof
for the special case of Theorem~\ref{satz:theo}
where $S$ consists of more than four points
in convex position, not three of them on a line.
The proof is based on a convergence argument similar in spirit to
our Lemma~\ref{lem:seq} and on bounds to distance ratios
between the converging points.

The rest of this paper is organized as follows.
In Section~\ref{sec:definitions} we list the exceptional
configurations
and we show basic properties of~$K(S)$.
In Section~\ref{sec:special} we study a special case of the problem.
In Section~\ref{sec:densekernel} we use arguments from projective
geometry to prove for any non-exceptional configuration~$S$ that
there exists a triangle~$T$ such that $\Sinfty$ is dense in~$T$.
Finally, in the last section we prove that in this case $\Sinfty$ is dense
in~$K(S)$.


\section{Exceptional configurations and the candidate}\label{sec:definitions}

\bigskip

Here, we list all cases of dilation-free graphs.
They can also be found at Eppstein's Geometry Junkyard~\cite{Junk}.
It can be proven by case analysis that these are all possibilities.
The {\em exceptional configurations} are the subsets of
the vertices of such graphs.
Most exceptional configurations are the
whole vertex set of a corresponding dilation-free graph.
Only the special case described in the footnote
in~($iii$) yields an
exceptional configuration which is a
proper subset.

\bigskip
\noindent
\setlength{\unitlength}{1mm}
\begin{picture}(130,6)

  \put(75,0){ \includegraphics{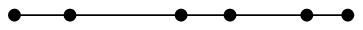}}
  \put(0,0)
  {
    \begin{minipage}[t]{60\unitlength}
    ($i$) $n$ points on a line
    \end{minipage}
  }
\end{picture}
\begin{picture}(130,30)
  \put(75,0){ \includegraphics{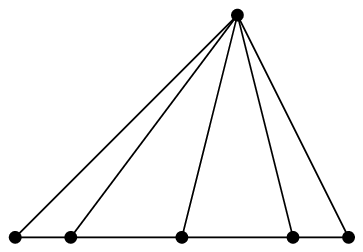}}
  \put(0,21)
  {
    \begin{minipage}[t]{60\unitlength}
      ($ii$) $n-1$ points on a line, one point not on this line
    \end{minipage}
  }
\end{picture}
\begin{picture}(130,30)
  \put(75,0){ \includegraphics{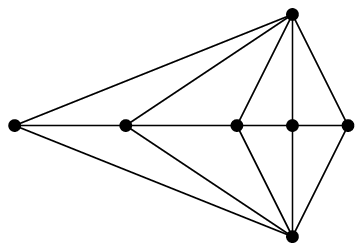}}
  \put(0,21)
  {
    \begin{minipage}[t]{60\unitlength}
      ($iii$) $n-2$ points on a line, two points on opposite sides of
      this line%
      \footnotemark
    \end{minipage}
  }
\end{picture}
\begin{picture}(130,30)
  \put(75,0){ \includegraphics{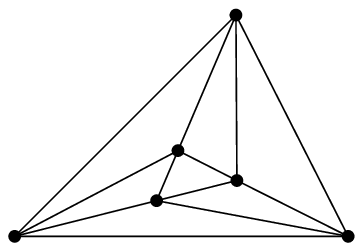}}
  \put(0,21)
  {
    \begin{minipage}[t]{60\unitlength}
      ($iv$) a triangle (i.e.\ three points) ne\-sted in the interior
      of another triangle.  Every pair of two inner points is
      collinear with one outer point.
    \end{minipage}
  }
\end{picture}
\footnotetext{Let $p_1$ and $p_2$ be the two points on opposite sides of
      the line, and let $p_3,\ldots,p_n$ be the other points.
      If the segment $p_1p_2$ intersects with
      the convex hull $\mbox{ch}(\{p_3,\ldots,p_n\})$,
      the intersection point must be a vertex of the dilation-free graph.
      However, it does not have to be part of the corresponding
      exceptional configuration.}

\bigskip
The main statement we want to prove in this paper is
that if~$S$ is not an exceptional configuration,
the set~$S_\infty$ is dense in the candidate~$K(S)$.
As a first step, we list some basic properties of~$K(S)$
in the following lemma.
The most important one is that there
cannot be any intersection point outside the candidate.

\begin{lemma}\label{lem:propertiescandidate}
  \begin{enumerate}
    \item[]
    \item\label{item:CandidatePolygon}
      If the candidate $K(S)$ is not empty, it is a convex polygon.
    \item\label{item:CandidateVertices}
      Every vertex of $K(S)$ belongs to $S_1 \cup S_2$.
    \item\label{item:IntersectionsInCandidate}
      Every intersection point lies inside of the candidate,
      that is $\Sinfty \setminus S \subseteq K(S)$.
  \end{enumerate}
\end{lemma}

\begin{proof}
  \ref{item:CandidatePolygon}.
  The convexity follows immediately from the definition
  because $K(S)$ is the intersection of (convex) closed half planes.
  The definition also implies that~$K(S)$ is
  a subset of the convex hull~$\ch(S)$.
  It is bounded.
  Furthermore, it is easy to see that in the definition of the candidate
  it suffices to consider only those finitely many closed half planes~$H$
  which have at least two points from~$S$ on their boundary.
  Hence, $K(S)$ is a polygon.

  \begin{figure}[htbp]
    \begin{center}
      \includegraphics{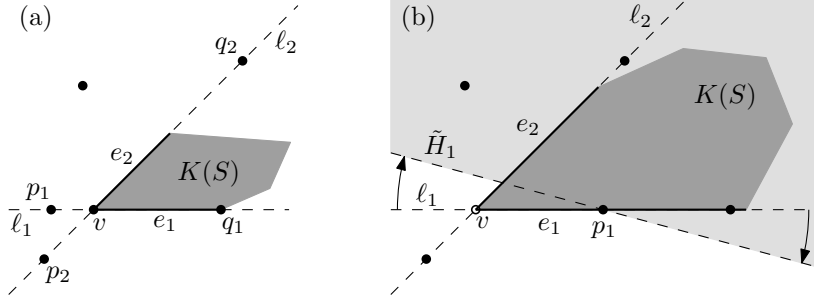}
      \caption{\label{fig:PropertiesCandidate} Every vertex~$v$
        of the candidate~$K(S)$
        has to be an element of $S_1 \cup S_2$.}
    \end{center}
  \end{figure}
  \noindent
  \ref{item:CandidateVertices}.
  Consider Figure~\ref{fig:PropertiesCandidate}.
  Let $e_1$ and $e_2$ be two neighbor edges of the polygon $K(S)$
  meeting in the vertex~$v$,
  and let $\ell_1$ and $\ell_2$ be the lines through these edges.
  By the arguments above there are at least two points from $S$
  on each of the lines.
  If there are two points $p_1, q_1 \in S$ on $\ell_1$
  on opposite sides of~$\ell_2$,
  and there are two points $p_2, q_2 \in S$ on $\ell_2$
  on opposite sides of~$\ell_1$,
  then, as shown in Figure~\ref{fig:PropertiesCandidate}.a,
  the vertex~$v$ belongs to $S_2$ since it is the
  crossing point of $p_1q_1$ and $p_2q_2$.

  Otherwise all points of $S \cap \ell_1$ lie on
  one side of~$\ell_2$ or all points of~$S \cap \ell_2$ lie
  on one side of~$\ell_1$.
  We assume the first situation,
  which is depicted in Figure~\ref{fig:PropertiesCandidate}.b.
  Let $p_1$ be the point in~$S \cap \ell_1$ which is closest to~$v$.
  Then, if $v$ is not a point of $S=S_1$,
  we can rotate the half plane $H_1$ belonging to $e_1$
  around $p_1$ such that the turned closed half
  plane~$\tilde{H}_1$ still contains all but one point of $S$
  and it does not contain $v \in K(S)$.
  This is a contradiction because by definition~$K(S)$
  would have to be contained in~$\tilde{H}_1$.

  \noindent
  \ref{item:IntersectionsInCandidate}.
  Let $H$ be a closed half plane containing $S\setminus\{p\}$.
  Then, clearly the crossing points
  of the next generation $S_2 \setminus S$
  are contained in $H$,
  since all the line segments connecting points from $S$ are either fully
  contained in $H$ or they meet in $p$.
  Applying this argument inductively
  shows that $\Sinfty \setminus S \subset H$.
\end{proof}


\section{A useful special case}\label{sec:special}

\begin{figure}[h]
  \begin{center}
  \includegraphics{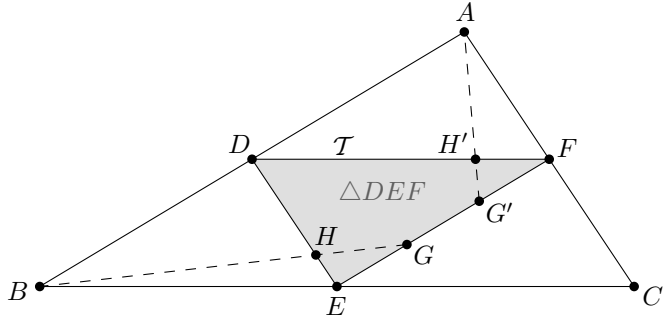}
  \caption{\label{fig:specialcase}
    Proving that the intersection points
    are dense in the triangle~$\triangle DEF$.}
  \end{center}
\end{figure}

In this whole section we consider a starting configuration~$S=\{A,B,C,D,E,F\}$
as shown in Figure~\ref{fig:specialcase}.
The three points~$A$, $B$, $C$ are not collinear,
and the three remaining points~$D$, $E$, $F$ are the midpoints
of the segments $\seg{A}{B}$, $\seg{B}{C}$ and $\seg{A}{C}$ respectively.
We want to prove the main statement for this special case.
This means, we want to show that for such a starting configuration
the set of intersection points~$\Sinfty$ is dense
in the triangle~$\triangle DEF$.

In order to do so,
it is sufficient to show that $\Sinfty$ is dense
on the sides of this triangle
because in this case the intersections of the segments connecting
the points on the sides are dense in the interior of the triangle.

Actually, for the proof we consider only a particular subset
$\Saltinfty \subseteq \Sinfty$.
Let~$\triangleboundary := \seg{D}{E} \cup \seg{E}{F} \cup \seg{D}{F}$
denote the boundary of the triangle~$\triangle{DEF}$.
Now, in each step of the iteration
we consider only segments
which connect points on~$\triangleboundary$
with a vertex from~$\{A,B,C\}$,
and the intersections of such segments with~$\triangleboundary$.
More formally we define
\begin{eqnarray*}
  &&
  \Salt_0
  :=
  \{D,E,F\},\;\;\;\;
  \forall i \in \N_0:\;\;  \tilde{L}_i
  :=
  \big\{
    pq
    \; \big| \;
    p \in \Salt_i,\;
    q \in \{A,B,C\}
  \big\}
  ,\\
  &&
  \Salt_{i+1}
  :=
  \big\{
    x
    \;\big|\;
    x \in l\cap \triangleboundary \mbox{ where }
    l \in \tilde{L}_i
  \big\},\;\;\;\;
  \mbox{and}\;\;\;
  \Saltinfty
  :=
  \bigcup_{i=1}^\infty \Salt_i.
\end{eqnarray*}
We call
a fraction~$k/l \in [0,1]$,
$k, l \in \N_0$, $\gcd(k,l)=1$,
\emph{constructible}
if there exists a point $H \in \Saltinfty$ on $\seg{D}{E}$
such that $|EH|/|DE| = k/l$.

The following observation, depicted in Figure~\ref{fig:specialcase},
contains the main step for deriving a new constructible number from
an old one.
\begin{lemma}\label{lem:constructiblestep}
  Let $G$ be a point on the segment~$\seg{E}{F}$
  such that $|EG|/|EF| = k/l$, 
  $k,l \in \N_0$.
  And let $H$ be the intersection
  of $\seg{B}{G}$ and $\seg{D}{E}$.
  Then, we have
  $|EH|/|DE| = k/(k+l)$.
\end{lemma}
\begin{proof}
  In the considered situation
  the segment~$\seg{B}{D}$ is parallel to $\seg{E}{F}$,
  which we denote by
  $\seg{B}{D} \parallel \seg{E}{F}$.
  Therefore, 
  the triangles $\triangle{EGH}$ and $\triangle{DBH}$
  are similar.
  This implies
  \begin{equation}\label{equ:fractionsintriangle}
    \frac{|DH|}{|EH|}
    =
    \frac{|BD|}{|EG|}
    =
    \frac{|EF|}{|EG|}
    =
    \frac{l}{k}
  \end{equation}
  The second equality holds because $DBEF$ is a parallelogram.
  We get
  \[
    \frac{|DE|}{|EH|}
    =
    \frac{|DH|+|HE|}{|EH|}
    \stackrel{\mbox{\scriptsize (\ref{equ:fractionsintriangle})}}{=}
    \frac{l}{k}+1
    =
    \frac{l+k}{k}.
  \]
\end{proof}

\noindent
We can prove this lemma analogously
for all other combinations of the sides of
the triangle~${\cal T}$ instead of $EF$ and $DE$.
In the next step we use this kind of symmetry to prove
that if we do not use $D$ and~$E$ in the definition
of constructible numbers
but a different pair of vertices of~$\triangle{DEF}$,
the set of constructible numbers remains the same.
Furthermore, we get two simple
rules for deriving these numbers.
Everything is summarized in the following
lemma.

\begin{lemma}\label{lem:constructible}
  \begin{enumerate}
    \item[]
    \item\label{lemsub:constructiblesymmetry}
      Let~$P, Q \in \{D,E,F\}$ be
      distinct vertices of~$\triangle{DEF}$.
      Then, $k/l$ is constructible if and only if
      there exists a point~$H \in \Saltinfty$ on
      $\seg{P}{Q}$ such that $|PH|/|PQ|=k/l$.
    \item\label{lemsub:constructiblemirror}
      The fraction~$k/l$ is constructible
      if and only if $1-k/l$ is constructible.
    \item\label{lemsub:constructiblestep}
      If~$k/l$ is constructible,
      then, also $k/(k+l)$ is constructible.
  \end{enumerate}
\end{lemma}

\begin{proof}
  \noindent
  \ref{lemsub:constructiblesymmetry}.~This
  is a rather obvious implication from the fact
  that we can prove Lemma~\ref{lem:constructiblestep}
  also for all other possible combinations of
  the sides of the triangle~$\triangle{DEF}$,
  not only~$\seg{E}{F}$ and $\seg{D}{E}$.

  Formally, one can prove by induction on~$i$
  that for every~$i \in \N_0$
  and every~$k/l \in [0,1]$,
  $k,l \in \N_0$, $\gcd(k,l)=1$, the following two statements
  are equivalent:
  \begin{itemize}
    \item[(a)] There exists a point~$H \in \Salt_i$ on $\seg{D}{E}$
      satisfying $|EH|/|DE|=k/l$.
    \item[(b)] For every pair of distinct vertices $P, Q \in \{D,E,F\}$,
      there exists a point~$H \in \Salt_i$ on $\seg{P}{Q}$
      satisfying $|PH|/|PQ|=k/l$.
  \end{itemize}
  The induction base~$i=0$ holds, because in this case
  (a) and (b) are both equivalent to $k/l \in \{0,1\}$.
  In the induction step the implication (b)~$\Rightarrow$~(a)
  is trivial.

  Now suppose the equivalence of~(a) and~(b) is already shown
  for~$i$.
  We have to prove ``(a)~$\Rightarrow$~(b)'' for~$i+1$.
  Assume (a) holds for $i+1$ and a given $k/l$.
  By definition of~$\Salt_{i+1}$
  there exists a point~$G \in \Salt_{i}$
  on~$\seg{E}{F}$ or on~$\seg{D}{E}$
  such that $H = \seg{B}{G} \cap \seg{D}{E}$.
  Here, we consider only the case where $G \in \seg{E}{F}$,
  $P=F$ and $Q=D$.
  All the other cases can be proved similarly.
  By Lemma~\ref{lem:constructiblestep}
  we know that~$|EG|/|EF| = k/(l-k)$.
  Consider Figure~\ref{fig:specialcase}.
  The induction hypothesis implies that there also exists
  a point $G' \in \Salt_i$ on $\seg{E}{F}$ satisfying $|FG'|/|EF|=k/(l-k)$.
  Now we apply a variant of Lemma~\ref{lem:constructiblestep}
  where we replace $G$ by $G'$, $H$ by $H'$, $E$ by $F$, 
  $F$ by $E$, $B$ by $A$, and~$D$ remains the same.
  This shows that
  the point~$H' := \seg{A}{G'} \cap \seg{D}{F} \in~\Salt_{i+1}$
  satisfies~$|FH'|/|DF| = k/(k+(l-k)) = k/l$.

  \noindent
  \ref{lemsub:constructiblemirror}.~This is an immediate
  consequence of the first statement.
  If we choose~$P:=D$ and $Q:=E$,
  it says that there exists a point~$H \in \Saltinfty$
  on~$\seg{D}{E}$ satisfying~$|DH|/|DE|=k/l$,
  hence $|DH| = |DE|(k/l)$.
  This implies
  \[
    |EH| = |DE|-|DH| = |DE|(1-k/l).
  \]

  \noindent
  \ref{lemsub:constructiblestep}.~Let~$k/l$ be constructible.
  Then,
  by choosing~$P=E$ and $Q=F$ in~\ref{lemsub:constructiblesymmetry}.,
  we derive
  that there exists a point~$G \in \Saltinfty$
  on $\seg{E}{F}$ such that $|EG|/|EF| = k/l$.
  We apply Lemma~\ref{lem:constructiblestep}.
  The resulting point $H := \seg{B}{G} \cap \seg{D}{E}$
  belongs to~$\Saltinfty$ and satisfies
  $|EH|/|DE| = k/(k+l)$.
  Hence,
  $k/(k+l)$ is constructible.
\end{proof}

\noindent
Now we can prove the main result of this section.
\begin{lemma}\label{lem:special}
  Let the starting configuration $S=\{A,B,C,D,E,F\}$ be as
  described in the beginning of this section,
  cf.\ Figure~\ref{fig:specialcase}.
  Then $\Sinfty$ is dense in $K(S)=\triangle{DEF}$.
\end{lemma}

\begin{proof}
  We will prove that every number in $[0,1] \cap \Q$
  is constructible.
  Of course, this is a dense set in~$[0,1]$.
  Hence, by the definition of constructible numbers
  and by Lemma~\ref{lem:constructible}.\ref{lemsub:constructiblesymmetry},
  this implies that~$\Saltinfty$ is dense in~$\triangleboundary$.
  Therefore, $\Sinfty$ is dense in~$\triangleboundary$,
  and thereby $\Sinfty$ is dense in~$\triangle{DEF}$.

  \sloppy Now, we prove that
  every $k/l \in [0,1]$, $k,l \in \N_0$, $\gcd(k,l)=1$,
  is constructible.
  We use induction on $n$, the sum of the numerator and the denominator
  of~$k/l$.

  If~$n=1$, we have $k=0$ and $l=1$. Therefore, $M=E$ proves the claim;
  the fraction~$0$ is constructible.
%
  Now, suppose the statement is true for the cases $k+l\leq n$.
  Consider an arbitrary fraction~$k/l \in (0,1]$
  where $k+l=n+1$ and $\gcd(k,l)=1$.

  \noindent
  {\em Case 1}, $0 < k/l \leq 1/2$:
  We define $k':= k$ and $l':=l-k$.
  We then have $k', l' \in \N_0$,
  $\gcd(k',l')=1$, $k'/l' \in (0,1]$ and $0 \leq k'+l' < k+l$.
  By induction hypothesis $k'/l'$ is constructible.
  Lemma~\ref{lem:constructible}.\ref{lemsub:constructiblestep}
  implies that $k'/(k'+l')=k/l$ is constructible,
  too.

  \noindent
  {\em Case 2}, $1/2 < k/l \leq 1$:
  We define $k':= l-k$ and $l':=l$.
  Then we have $k', l' \in \N_0$, $\gcd(k',l')=1$,
  $k'/l' \in [0,1)$ and $k'+l' < k+l$.
  Thus by induction hypothesis $k'/l'$ is constructible.
  By Lemma~\ref{lem:constructible}.\ref{lemsub:constructiblemirror}
  we can also construct $k/l = 1 - k'/l'$.
\end{proof}


\section{Density in a triangle}\label{sec:densekernel}

In this section we want to prove that whenever
$S$ is not an exceptional configuration,
the set of intersection points~$\Sinfty$ is dense
in a triangle.
We start with a special case.

\begin{figure}[h]%
  \begin{center}%
    \includegraphics{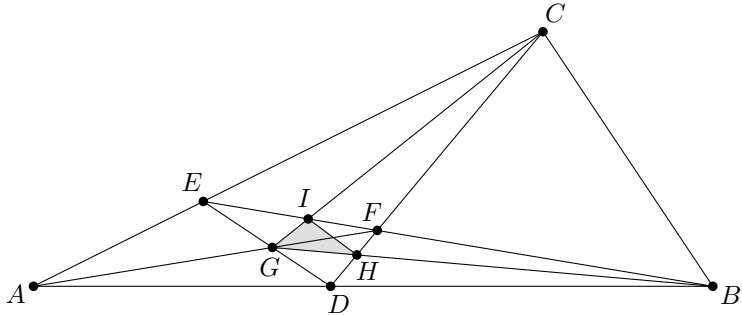}%
    \caption{\label{fig:denseintriangle}%
      $\Sinfty$ is dense in $\triangle{GHI}$.}%
  \end{center}%
\end{figure}
%
\begin{lemma}\label{lem:main} (Main Lemma)
  Let the starting configuration $S=\{A,B,C,D,E\}$ be as
  shown in Figure~\ref{fig:denseintriangle}.
  The three points $A$, $B$, $C$ are not collinear.
  The point $D$ lies on the segment~$\seg{A}{B}$,
  and $E$ is located on $\seg{A}{C}$.
  Then $\Sinfty$
  is dense in a triangle. 
\end{lemma}
\begin{figure}[hbt]%
  \begin{center}%
    \includegraphics{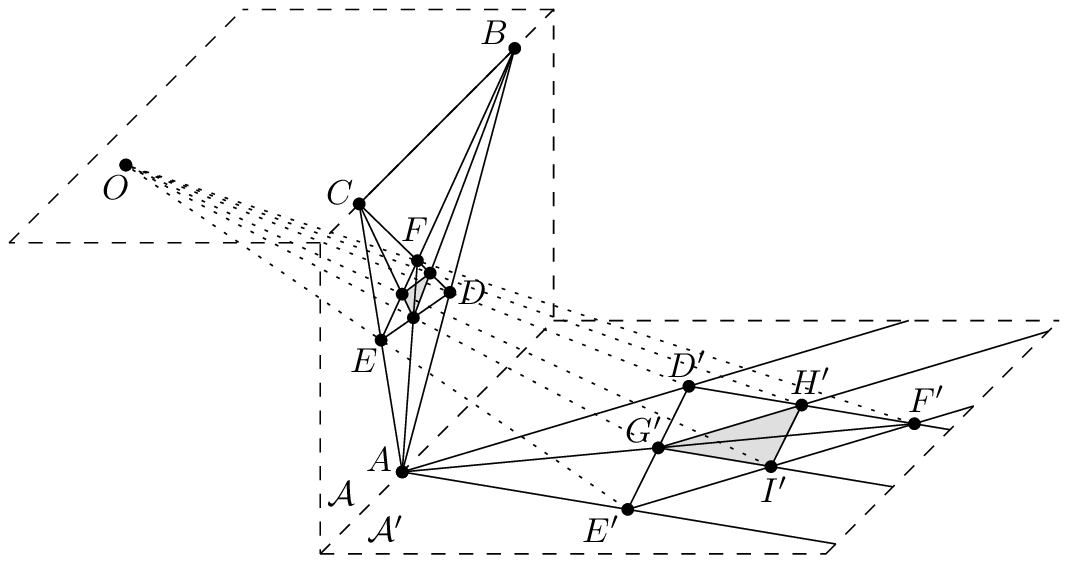}%
    \caption{\label{fig:projectiontriangle}
      A useful projection of $\triangle ABC$.}%
  \end{center}%
\end{figure}

\begin{proof}
  We use the projection shown in Figure~\ref{fig:projectiontriangle}.
  The given two-dimensional geometric situation is
  represented in an affine plane~$\A$ of~$\R^3$
  which does not contain the origin~$O$.
  Additionally, we consider the affine plane~$\A'$ which
  is parallel to~$\seg{O}{B}$ and $\seg{O}{C}$,
  and which contains~$A$.

  We can define a projection~$\pi: \triangle{ABC} \setminus \seg{B}{C} \to \A'$
  by mapping each point~$P \in \triangle{ABC} \setminus \seg{B}{C}$
  to the one point in~$\A'$
  which is hit by the line through~$O$ and $P$.
  For every such point~$P$ we denote $\pi(P)$ by $P'$.

  Basic arguments from projective geometry show that
  the projection~$\pi$ is a homeomorphism
  between~$\triangle{ABC} \setminus \seg{B}{C}$
  and $\pi(\triangle{ABC} \setminus \seg{BC})$,
  and that it maps line segments to line segments.
  Furthermore, the projections of two segments
  which meet in~$B$
  are parallel in~$\A'$,
  and the same holds for two segments which meet in~$C$.
  For an introduction to projective geometry
  see for instance Bourbaki~\cite{Bourb}.

  In~$\A'$ we have $\seg{A'}{D'} \parallel \seg{E'}{F'}$,
  because the lines through~$\seg{A}{D}$ and~$\seg{E}{F}$
  intersect in~$B$,
  and analogously $\seg{A'}{E'} \parallel \seg{D'}{F'}$,
  because the corresponding lines intersect in~$C$.
  Hence $G'$ is the midpoint of
  $\seg{D'}{E'}$, $I'$ is the midpoint of $\seg{D'}{F'}$,
  and $H'$ is the midpoint of~$\seg{E'}{F'}$.
  By Lemma~\ref{lem:special},
  for the starting configuration~$S' := \{D',E',F',G',H',I'\}$
  the iterated crossing points~$\Sinfty'$ are dense
  in the triangle~$\triangle{G'H'I'}$.
  Because~$\pi$ is a homeomorphism,
  this shows that~$\Sinfty$ is dense in~$\triangle{GHI}$.
\end{proof}

\begin{corollary}\label{cor:neck}
  Let the starting configuration~$S$
  be a set of $n>4$ points in the plane in convex position,
  no three of them on a line.
  Then $\Sinfty$ is dense in a triangle.
\end{corollary}
\begin{figure}[hbt]
  \begin{center}%
    \includegraphics{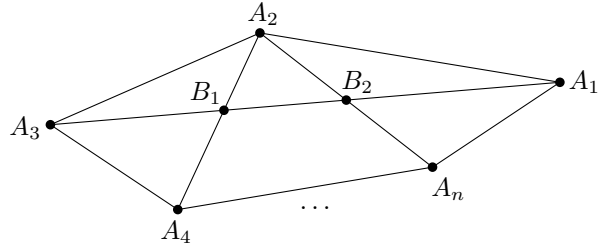}%
    \caption{\label{fig:n-gon} Proving the density in a triangle for
      $n>4$ points in convex position.}%
  \end{center}%
\end{figure}%

\begin{proof}
  Consider Figure~\ref{fig:n-gon}.
  We label the starting points $A_1,A_2,\ldots, A_n$
  counterclockwise.
  We define $B_1:= \seg{A_1}{A_{3}}\cap \seg{A_{2}}{A_{4}}$
  and $B_2:= \seg{A_1}{A_{3}}\cap \seg{A_{2}}{A_{n}}$.
  Now we can apply the main lemma, Lemma~\ref{lem:main},
  to the triangle $\triangle{A_{2}A_4A_n}$
  with $B_1$ and $B_2$ on different sides of this triangle.
  Thus $\Sinfty$ is dense in a triangle.
\end{proof}

\begin{lemma}\label{lem:densekernel}
  For any non-exceptional configuration,
  there exists a triangle in which $\Sinfty$ is dense.
\end{lemma}

\begin{proof}
  If the convex hull of $S_1$,
  which we denote by~$\ch(S_1)$,
  has more than four vertices,
  the claim is valid by Corollary~\ref{cor:neck}.
  If the convex hull has three vertices $A$, $B$ and $C$,
  we consider the following cases:

  \noindent
  {\bf Case 1}:
  There exist two additional starting points on
  different sides of the triangle. In this case the claim is also
  valid by the main lemma.
  \begin{figure}[h]%
    \begin{center}%
      \includegraphics{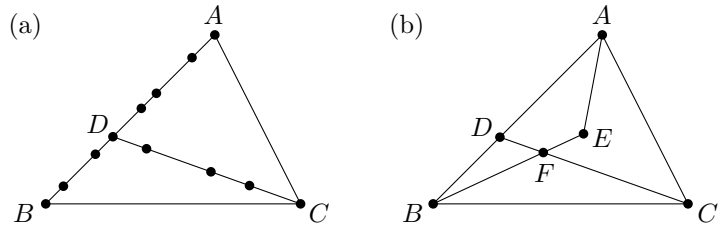}%
      \caption{\label{fig:case2} Case 2: $|\ch(S)|=3$ and there is a
        starting point on one side of the triangle.}%
    \end{center}%
  \end{figure}

  \smallskip\noindent
  {\bf Case 2}:
  There exists a starting point $D$ on one
  side of triangle $\triangle{ABC}$, say on $\seg{A}{B}$.
  The configuration contains at least one additional point,
  since it is not exceptional,
  but none on $\seg{A}{C}$ nor on $\seg{B}{C}$.
  If the rest of the points $S\setminus \{A,B,C,D\}$
  are all located on the line segment~$\seg{A}{B}$,
  or all on~$\seg{C}{D}$, then the
  configuration belongs to the exceptional case~(ii) or~(iii).
  Hence we have one of the following cases,
  as shown in Figure~\ref{fig:case2}.

  \smallskip\noindent
  (a) All of the additional points $S\setminus \{A,B,C,D\}$
  lie on $\seg{A}{B}$ and $\seg{C}{D}$, and on both line
  segments exists at least one additional point.
  In this case we can apply the main lemma to at least one of the triangles
  $\triangle{ACD}$ or~$\triangle{BCD}$. 

  \smallskip\noindent
  (b) There exists at least one additional point~$E$ in the interior
  of the triangle~$\triangle{ABC}$, but not on $\seg{C}{D}$.
  In this case, $E$ is either in
  the interior of triangle~$\triangle{ACD}$
  or in the interior of~$\triangle{BCD}$.
  Hence one of the line segments $\seg{B}{E}$ or $\seg{A}{E}$
  intersects~$\seg{C}{D}$ in a point~$F$.
  Now we can apply the main lemma to triangle~$\triangle{ABE}$
  with $D$ and $F$ on different sides. 
  \begin{figure}[hbt]%
    \begin{center}%
      \includegraphics{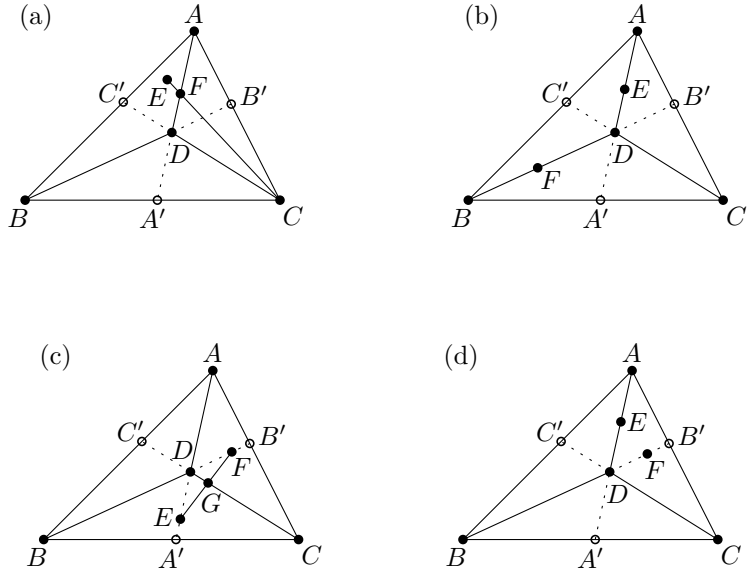}
      \caption{\label{fig:case3} Case 3: $|\ch(S)|=3$
        and there is no additional point on the boundary of~$\triangle{ABC}$.}%
    \end{center}%
  \end{figure}

  \smallskip\noindent
  {\bf Case 3}:
  There is no additional point on the boundary of triangle~$\triangle{ABC}$.
  Hence there exists at least one point $D$
  in the interior of $\triangle{ABC}$.
  The configuration requires at least one more point not to be an
  exceptional configuration.
  Let $A'$ be the intersection point of $BC$ and the line through $A$
  and $D$.
  We define $B'$ and $C'$ analogously.
  Note that these points
  are not necessarily in $\Sinfty$!

  If the rest of the points are all on $\seg{A}{A'}$,
  all on $\seg{B}{B'}$, or all on~$\seg{C}{C'}$,
  we obtain the exceptional configuration~(iii).
  Otherwise we distinguish the following four cases,
  as shown in Figure~\ref{fig:case3}.

  \smallskip\noindent
  (a)~There exists a point~$E$ in the interior of triangle~$\triangle{ABC}$
  which is not located on $\seg{A}{A'}$, $\seg{B}{B'}$ or $\seg{C}{C'}$.
  Without loss of generality
  it is in the interior of triangle~$\triangle{ABD}$.
  Thus $\seg{C}{E}$ intersects~$\seg{A}{D}$ or~$\seg{B}{D}$ in a point~$F$;
  we may assume that $F$ lies on~$\seg{A}{D}$.
  Now we have the same situation as in Case 2.b for the
  triangle $\triangle{ABD}$ with point~$F$ on a side and point~$E$ in
  the interior of~$\triangle{ABD}$.
  The points $B$, $E$ and~$F$ cannot be collinear
  because $C$, $E$ and~$F$ are collinear by construction
  and neither~$E$ nor~$F$ lies on~$\seg{B}{C}$.

  \smallskip\noindent
  (b)~At least two of the segments $\seg{A}{D}$,
  $\seg{B}{D}$ and $\seg{C}{D}$ contain an additional point.
  We may assume that there is a point~$E$ on~$\seg{A}{D}$
  and $F$ on~$\seg{B}{D}$.
  Then we can apply the main lemma to triangle $\triangle{ABD}$.

  \smallskip\noindent
  (c)~At least two of the segments $\seg{A'}{D}$, $\seg{B'}{D}$ and
  $\seg{C'}{D}$ contain an additional point.
  We may assume that there is a point~$E$ on~$\seg{A'}{D}$
  and a point~$F$ on~$\seg{B'}{D}$.
  Then $E$ lies in the interior of~$\triangle{BCD}$,
  and $F$ lies in the interior of~$\triangle{ACD}$.
  The segment~$\seg{E}{F}$ intersects~$CD$ in a point~$G$,
  since $A'DB'C$ is convex.
  Now we have the case 2.b for triangle $\triangle{ACD}$
  with $G$ on one side and $F$ in the interior.
  If the points $A$, $F$ and $G$ were collinear,
  all of them had to be located on~$\seg{A}{A'}$,
  since $E$, $F$ and $G$ are collinear,
  and $A$ and $E$ lie on $\seg{A}{A'}$.
  This contradicts $F$ being contained
  in the interior of~$\triangle{ACD}$.

  \smallskip\noindent
  (d)~One of the segments $\seg{A}{D}$, $\seg{B}{D}$ or $\seg{C}{D}$
  contains an additional point~$E$,
  and one of the segments $\seg{A'}{D}$, $\seg{B'}{D}$ or $\seg{C'}{D}$
  contains an additional point~$F$.
  We may assume $E\in \seg{A}{D}$.
  Now, if we had only additional points on~$\seg{A'}{D}$,
  this would be exceptional case~(iii).
  Hence we may assume~$F \in \seg{B'}{D}$.
  The other cases can be treated analogously.
  All the starting points in $S \setminus \{A,B,C,D\}$
  lie on~$\seg{A}{D}$ or~$\seg{B'}{D}$
  since otherwise we had case~3.(a), (b) or (c).
  If we have only one additional starting point~$E$
  on~$\seg{A}{D}$ and only one point~$F$ on~$\seg{B'}{D}$,
  and $C$, $E$ and $F$ are collinear,
  then we have exceptional configuration~(iv).
  Otherwise, there exist starting points~$E$ on~$\seg{A}{D}$
  and~$F$ on~$\seg{B'}{D}$ such that $C$, $E$ and $F$ are not
  collinear.
  We have case~2.(b) for triangle~$\triangle{ACD}$.

  \begin{figure}[h]%
    \begin{center}%
      \includegraphics{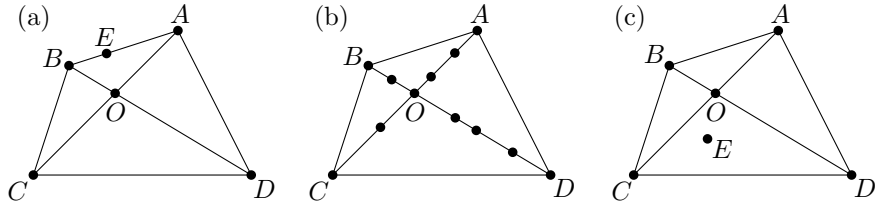}
      \caption{\label{fig:FourVerticesOnCH} Subcases, if the convex hull of a
        non-exceptional starting configuration has four vertices.}%
    \end{center}%
  \end{figure}%
  In the remaining case the convex hull has four vertices,
  say $A$, $B$, $C$ and $D$ in counterclockwise order.
  Let $O$ denote the intersection point of~$\seg{A}{C}$ and~$\seg{B}{D}$.
  If the additional points of~$S\setminus\{A,B,C,D\}$
  lie all on~$\seg{A}{C}$ or all on~$\seg{B}{D}$,
  we have exceptional case~(iii).
  Otherwise we distinguish three cases,
  as shown in Figure~\ref{fig:FourVerticesOnCH}.

  \smallskip\noindent
  (a)~If there exists an additional point $E$ on one of the sides
  of the quadrilateral~$ABCD$,
  say on~$\seg{A}{B}$,
  we can apply the main lemma to triangle~$\triangle{ABC}$
  with $E$ on one side and $O$ on the other side.

  \smallskip\noindent
  (b)~If the rest of the points are all located on $AC$ and $BD$ but
  not only on one of them, then we can also use the main lemma for at
  least one of the triangles $\triangle{ABO}, \triangle{ADO}, \triangle{BCO}$ or $\triangle{CDO}$.

  \smallskip\noindent
  (c)~If there is a point $E$ in the interior of the quadrilateral but
  not on~$AC$ or~$BD$, then $E$ lies in one of the triangles
  $\triangle{ABC}$ or $\triangle{ACD}$. This is case 2.b for the
  triangle containing $E$.
\end{proof}

\section{Density in the candidate}\label{sec:theo}

\begin{figure}[hbt]%
  \begin{center}%
    \includegraphics{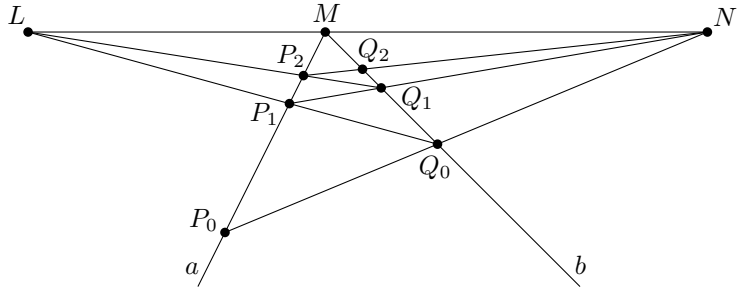}%
    \caption{\label{fig:seq} The point sequences~$(P_i)_{i\in\N_0}$
      and~$(Q_i)_{i\in\N_0}$ converge to~$M$.}%
  \end{center}%
\end{figure}

We now know that, if~$S$ is not an exceptional configuration,
the intersection points~$\Sinfty$ are dense in a triangle.
We still want to show that they are dense in~$K(S)$.
To this end, we introduce some lemmata.
The situation of Lemma~\ref{lem:seq}
is shown in Figure~\ref{fig:seq}.
\begin{lemma}\label{lem:seq}
  Let $L$, $M$ and $N$ be three distinct points
  such that $M$ lies on the line segment~$\seg{L}{N}$,
  and let $a \mbox{ and }b$ be two rays emanating
  from $M$ on the same side of~$\seg{L}{M}$.
  Let $P_0$ be a point on the ray~$a$.
  We define a sequence of points inductively by
  $Q_i := \seg{P_i}{N}\cap b$,
  $P_{i+1} := \seg{Q_i}{L}\cap a$
  for every~$i \in \N_0$.
  Then the point sequence $(P_i)_{i\in\N_0}$
  converges to~$M$ on~$a$
  and the point sequence $(Q_i)_{i\in\N_0}$
  converges to~$M$ on~$b$.
\end{lemma}
\begin{figure}[hbt]%
  \begin{center}%
    \includegraphics{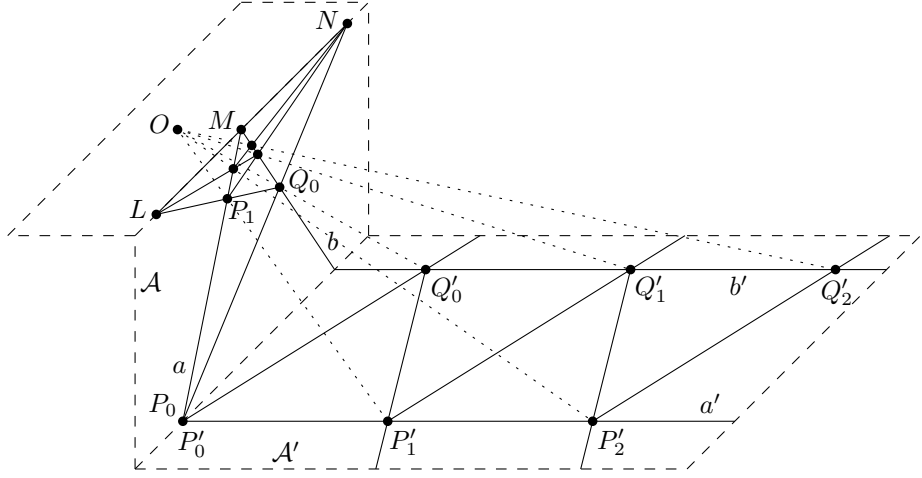}%
    \caption{\label{fig:seqproj} Proving Lemma~\ref{lem:seq}
      with a projection.}%
  \end{center}%
\end{figure}%
\begin{proof}
  One way to prove this simple result
  is to apply the same trick as in the proof of the main lemma.
  We use the projection shown in Figure~\ref{fig:seqproj}.
  The given two-dimensional geometric situation is
  represented in an affine plane~$\A$ of~$\R^3$
  which does not contain the origin~$O$.
  Additionally, we consider an affine plane~$\A'$ which
  is parallel to~$\seg{O}{L}$ and $\seg{O}{N}$,
  and which contains~$P_0$.

  We define a
  projection~$\pi: \triangle{P_0LN} \setminus \seg{L}{N} \to \A'$
  by mapping each point~$P \in \triangle{P_0LN} \setminus \seg{L}{N}$
  to the one point in~$\A'$
  which is hit by the line through~$O$ and~$P$.
  Again for every such point~$P$ we denote $\pi(P)$ by $P'$.

  And basic arguments from projective geometry show that
  the projection~$\pi$ is injective and maps line segments to line segments.
  Furthermore, the projections of two segments
  which meet in~$M$ are parallel in~$\A'$,
  and the same holds for two segments which meet in~$L$ or~$N$.

  Because of this, we know that the projected rays~$a' := \pi(a)$
  and $b' := \pi(b)$ are parallel in $\mathcal{A}'$.
  With the same argument we get
  \begin{displaymath}
    \seg{P_0'}{Q_0'}
    \parallel
    \seg{P_1'}{Q_1'}
    \parallel
    \cdots
    \parallel
    \seg{P_i'}{Q_i'}
    \parallel
    \cdots
    \mbox{ and }
    \seg{Q_0'}{P_1'}
    \parallel
    \seg{Q_1'}{P_2'}
    \parallel
    \cdots
    \parallel
    \seg{Q_i'}{P_{i+1}'}
    \parallel
    \cdots.
  \end{displaymath}
  Hence $P_i'P_{i+1}'Q_{i+1}'Q_{i}'$ and $Q_i'Q_{i+1}'P_{i+2}'P_{i+1}'$
  are parallelograms for every $i\in \N_0$.
  This shows
  \begin{displaymath}
    |P_0'P_1'|
    =
    |Q_0'Q_1'|
    =
    |P_1'P_2'|
    =
    |Q_1'Q_2'|
    =
    \ldots .
  \end{displaymath}
  These equations imply that the sequences $(P_i')_{i\in\N_0}$
  and $(Q_i')_{i\in\N_0}$ diverge,
  they converge to the endpoint at infinity
  of $a'$, $b'$ respectively.
  This shows that $(P_i)_{i\in\N_0}$
  converges to~$M$,
  because otherwise it had to converge to a different
  point on~$a$,
  but then $(P_i')_{i\in\N_0}$ had to converge
  to the corresponding point on~$a'$.
  Analogously we derive that $(Q_i)_{i \in \N_0}$ converges to~$M$.
\end{proof}

\begin{figure}[hbt]%
  \begin{center}%
    \includegraphics{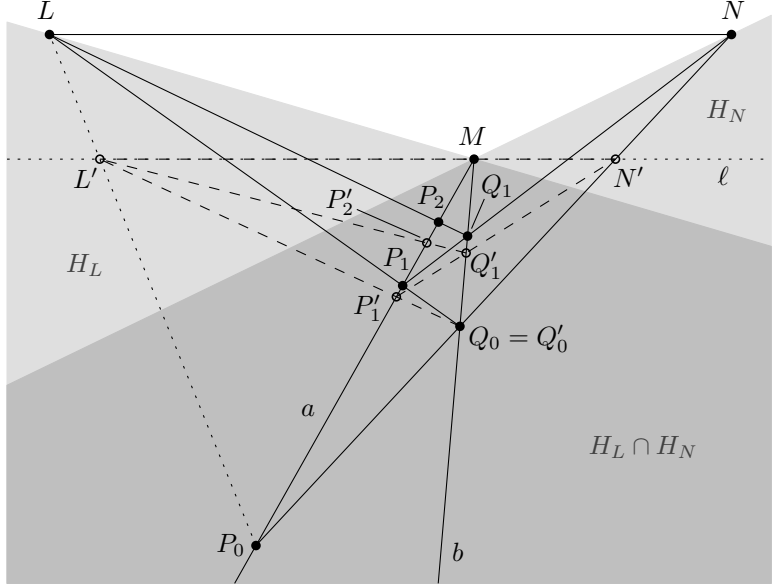}%
  \caption{\label{fig:seqcor} These point sequences~$(P_i)_{i\in\N_0}$
    and~$(Q_i)_{i\in\N_0}$ also converge to~$M$.}%
  \end{center}%
\end{figure}

\noindent
We can generalize this statement to the situation
shown in Figure~\ref{fig:seqcor}.
For two distinct points~$P$, $Q$, let~$\line(P,Q)$
denote the line which passes through both points.
\begin{corollary}\label{cor:seq}
  Let $L$ and $N$ be two points,
  and this time let $M$ be a point not on the segment~$\seg{L}{N}$
  and not collinear with~$L$ and $N$.
  Let $H_L$ be the half plane bounded by~$\line(L,M)$
  which does not contain~$N$,
  and similarly let~$H_N$ be the half plane bounded by~$\line(M,N)$
  which does not contain~$L$.
  Let $a$ and~$b$ be two rays emanating from~$M$ and lying in
  the interior of~$H_L \cap H_N$.
  If~$P_0$ is a point on the ray~$a$ and
  we define point sequences~$(P_i)_{i\in\N_0}$ and~$(Q_i)_{i\in\N_0}$
  analogously to Lemma~\ref{lem:seq},
  then both sequences converge to~$M$.
\end{corollary}

\begin{proof}
  Consider Figure~\ref{fig:seqcor}.
  Let~$\ell$ be the line through~$M$ which is parallel to $\line(L,N)$.
  Then we define~$L' := \ell \cap \seg{P_0L}$ and
  $N' := \ell \cap{P_0N}$.
  Now, because~$M$ lies on the segment~$\seg{L'}{N'}$,
  we can apply Lemma~\ref{lem:seq} to~$L'$, $M$, $N'$, $a$, $b$
  and $P_0':=P_0$ to construct two convergent sequences,
  $(P'_i)_{i\in\N_0}$ on~$a$ and $(Q'_i)_{i\in\N_0}$ on~$b$,
  which `push'~$(P_i)_{i\in\N_0}$ and $(Q_i)_{i\in\N_0}$ towards $M$,
  that is $|P_iM| \leq |P'_iM| \searrow 0$
  and $|Q_iM| \leq |Q'_iM| \searrow 0$.
\end{proof}

\begin{figure}[hbt]%
  \begin{center}%
    \includegraphics{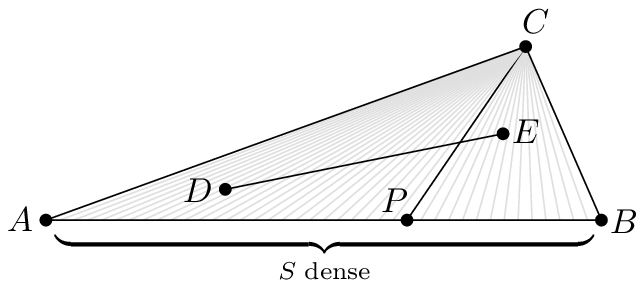}
    \caption{\label{fig:densityonsegment} If~$S$ is dense on~$\seg{A}{B}$,
      the intersections
      $\{\seg{PC}\cap\seg{DE}\;|\; P \in S\}$
      are dense on~$\seg{D}{E}$.}%
  \end{center}%
\end{figure}

\noindent
The last technical tool we need is the following fact,
depicted in Figure~\ref{fig:densityonsegment}.
We omit the simple proof.
\begin{lemma}\label{lem:densityonsegment}
  Let~$\triangle{ABC}$ be a non-degenerate triangle,
  i.e.\ $A$, $B$ and~$C$ are not collinear,
  and let~$D$ and $E$ be two points in~$\triangle{ABC}$
  which are not collinear with~$C$.
  If we connect every point of a point set~$S$,
  which is dense on~$\seg{A}{B}$,
  with~$C$,
  then the intersections of the resulting line segments
  with~$\seg{D}{E}$ are dense on~$\seg{D}{E}$.
\end{lemma}

\noindent
Finally, before we can prove the main result,
we introduce some additional useful notation.
\begin{figure}[hbt]%
  \begin{center}%
    \includegraphics{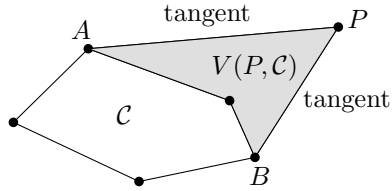}%
    \caption{\label{fig:vis} The visibility cone of~$P$
      with respect to $\mathcal{C}$, $V(P,\mathcal{C}$).}%
  \end{center}%
\end{figure}%
\begin{definition}
  Let $\mathcal{C}$ be a convex polygon and let $P$ be a point outside
  of $\mathcal{C}$.
  As before, for any set~$S \subset \R^2$,
  we write~$\ch(S)$ to denote the convex hull of~$S$.
  We define the {\em visibility cone} of~$P$ with respect to~$\mathcal{C}$
  by
  \begin{displaymath}
    V(P,\mathcal{C}):= \ch(\{P\}\cup\mathcal{C}) \setminus \mathcal{C}.
  \end{displaymath}
  We call the two edges of~$V(P,\mathcal{C})$ which are adjacent to~$P$
  the {\em tangents} of~$\mathcal{C}$ through~$P$.
\end{definition}

\smallskip\noindent
Now, we are ready to prove the main result of this article.
The following lemma contains the missing argument.
\begin{lemma}\label{lem:develop}
  Let the starting configuration $S\subset \R^2$
  be an arbitrary finite point set.
  Assume that~$\Sinfty$ is dense in a
  convex polygon~$\mathcal{C} \subset \ch(S)$ with nonempty interior.
  And let $M$ be a point from~$\Sinfty$
  such that $M \not\in \mathcal{C}$ and~$M$ is not a vertex of~$\ch(S)$.
  Then $\Sinfty$ is dense in~$\ch(\{M\} \cup \mathcal{C})$.
\end{lemma}
\begin{figure}[hbt]%
  \begin{center}%
    \includegraphics{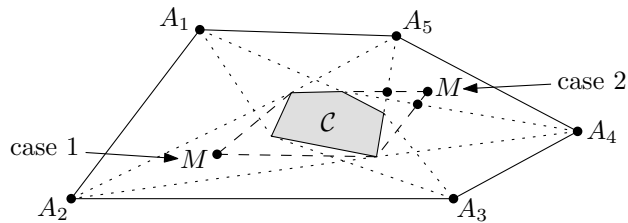}%
    \caption{\label{fig:extendingdensity} The two cases
      in the proof of Lemma~\ref{lem:develop}.}%
  \end{center}%
\end{figure}%

\begin{proof}
  Let $A_1,\ldots,A_n$ be the vertices of~$\ch(S)$
  in counterclockwise order.
  We distinguish the following two cases,
  as shown in Figure~\ref{fig:extendingdensity}.

  \smallskip\noindent
  {\em Case 1:}
  If the point~$P$ lies in the interior of
  a visibility cone~$V(A_i, \mathcal{C})$,
  we have $V(M,\mathcal{C}) \subset V(A_{i}, \mathcal{C})$.
  Thus by Lemma~\ref{lem:densityonsegment},
  $\Sinfty$ is dense on both tangents of~$\mathcal{C}$
  through~$P$.
  Clearly, together with the density in~$\mathcal{C}$,
  this implies that~$\Sinfty$ is dense in $V(M, \mathcal{C})$,
  and therefore it is dense in $\ch(\{M\} \cup \mathcal{C})$.

  If~$M$ lies on the boundary of a visibility cone~$V(A_i, \mathcal{C})$,
  Lemma~\ref{lem:densityonsegment} shows only that~$\Sinfty$
  is dense on the one tangent which is not collinear with~$A_i$.
  But in this special case this is enough to prove
  the density in~$\ch(\{M\} \cup \mathcal{C})$.

  \begin{figure}[hbt]%
    \begin{center}%
      \includegraphics{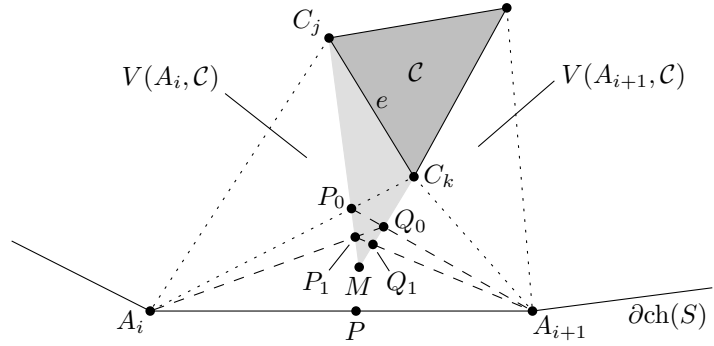}%
      \caption{\label{fig:extendingdensitycase2}
        If~$M$ does not lie inside a cone~$V(A_i,\mathcal{C})$,
        $\Sinfty$ is dense in~$V(M,\mathcal{C})$
        because of Lemma~\ref{lem:densityonsegment}
        and Corollary~\ref{cor:seq}.}%
    \end{center}%
  \end{figure}%
  \smallskip\noindent
  {\em Case 2:}
  Consider Figure~\ref{fig:extendingdensitycase2}.
  The visibility cone~$V(M,\mathcal{C})$ intersects
  at least one~$V(A_i,\mathcal{C})$.
  This can be seen by considering a point~$P$
  continuously moving along the boundary of~$\ch(S)$.
  At some point its visibility cone~$V(P,\mathcal{C})$
  contains~$M$.
  If $P$ is a vertex of~$\ch(S)$, we have case~1.
  Otherwise, if we define~$A_{n+1} := A_1$,
  the point~$P$ lies on a segment~$\seg{A_i}{A_{i+1}}$,
  $i \in \{1, \ldots, n\}$.

  The point~$M$ sees at least one whole edge~$e$ of~$\mathcal{C}$.
  Because of~$V(M,\mathcal{C}) \subseteq V(P,\mathcal{C})$,
  also~$P$ can see~$e$.
  The edge~$e$ is also visible from~$A_i$ or from~$A_{i+1}$.
  This can be proved by moving a point~$Q$ along~$\seg{A_i}{A_{i+1}}$.
  It cannot be that~$e$ is first invisible from~$Q$, then becomes visible,
  and gets invisible again.
  Because~$e$ is visible from~$P$, it therefore must also be visible from
  $A_i$ or~$A_{i+1}$.
  This proves that~$V(M,\mathcal{C})$ intersects $V(A_i,\mathcal{C})$
  or~$V(A_{i+1},\mathcal{C})$.

  We assume that it intersects~$V(A_{i},\mathcal{C})$.
  And let~$\seg{M}{C_j}$ and~$\seg{M}{C_k}$ be the tangents
  of~$\mathcal{C}$ through~$M$.
  Then, a part of one of the tangents, say of~$\seg{M}{C_j}$,
  lies inside~$V(A_{i},\mathcal{C})$.
  Let~$P_0$ be the intersection point of $\seg{M}{C_j}$
  and the boundary of~$V(A_i,\mathcal{C})$.
  Then, by Lemma~\ref{lem:densityonsegment} the set~$\Sinfty$
  is dense on~$\seg{P_0}{C_j}$.
  Let~$Q_0$ be the intersection of~$\seg{M}{C_k}$
  and $\seg{A_{i+1}}{P_0}$.
  Then, by Lemma~\ref{lem:densityonsegment} the set~$\Sinfty$ is dense on
  $\seg{Q_0}{C_k}$.
  We can use this argument inductively.
  The next step is to define~$P_1 := \seg{A_i}{Q_0} \cap \seg{M}{C_j}$.
  and~$Q_1 := \seg{A_{i+1}}{P_1} \cap \seg{M}{C_k}$.
  Lemma~\ref{lem:densityonsegment} shows that~$\Sinfty$
  is dense on~$\seg{P_1}{P_0}$, and that it is dense on~$\seg{Q_1}{Q_0}$.
  Corollary~\ref{cor:seq} proves that the sequences defined this way,
  $(P_i)_{i \in \N_0}$ and $(Q_i)_{i \in \N_0}$,
  both converge to~$M$.

  Hence,~$\Sinfty$ is dense on both tangents, $\seg{M}{C_j}$ and $\seg{M}{C_k}$.
  Because it is also dense in~$\mathcal{C}$,
  this shows that it is dense in~$\ch(M \cup \mathcal{C})$.
\end{proof}

In the end,
we only have to put all the pieces together to prove
that for every non-exceptional starting configuration~$S$
the set of intersection points~$S_{\infty}$ is dense in~$K(S)$.
\begin{proof} (of Theorem~\ref{satz:theo})
  If $S$ is not an exceptional configuration,
  by Lemma~\ref{lem:densekernel} we know
  that $\Sinfty$ is dense in a triangle $T$.
  This implies by
  Lemma~\ref{lem:propertiescandidate}.\ref{item:IntersectionsInCandidate}
  that $T \subseteq K(S)$.
  We also know by~Lemma~\ref{lem:propertiescandidate}
  that $K(S)$ is a convex polygon whose vertices belong to~$\Sinfty$.
  Let $P_1, P_2,\ldots, P_n$ be the vertices of~$K(S)$.
  Then, by using Lemma~\ref{lem:develop},
  we can show that~$\Sinfty$ is dense in~$\ch(T\cup \{P_1\})$.
  By applying Lemma~\ref{lem:develop} again,
  we get that~$\Sinfty$ is dense in~$\ch(T\cup\{P_1,P_2\})$.
  We repeat this argument $n$ times to prove that
  $\Sinfty$ is dense in~$\ch(T_\cup\{P_1,\ldots,P_n\})=K(S)$.
\end{proof}

\section{Acknowledgements}
We thank Dr. Behrang Noohi for his nice idea to use projective
geometry for the proof of the main lemma
and the anonymous referees for their very helpful and detailed
comments.


\begin{thebibliography}{10}

\bibitem{Frankl}
I.~B\'{a}r\'{a}ny, P.~Frankl, and H.~Maehara.
\newblock Reflecting a triangle in the plane.
\newblock {\em Graphs and Combinatorics}, 9:97--104, 1993.

\bibitem{Pach}
K.~Bezdek and J.~Pach.
\newblock A point set everywhere dense in the plane.
\newblock {\em Elemente der Mathematik}, 40(4):81--84, 1985.

\bibitem{Bourb}
N.~Bourbaki.
\newblock {\em General Topology, part 2}, chapter VI, {\S}3, pages 44--53.
\newblock Elements of Mathematics. Hermann, Paris, 1966.

\bibitem{Klein}
A.~Ebbers-Baumann, A.~Gr{\"u}ne, M.~Karpinski, R.~Klein, C.~Knauer, and
  A.~Lingas.
\newblock Embedding point sets into plane graphs of small dilation.
\newblock In {\em Algorithms and Computation: 16th International Symposium,
  ISAAC 2005}, volume 3827 of {\em Lecture Notes Comput. Sci.}, pages 5--16.
  Springer, December 2005.

\bibitem{Junk}
D.~Eppstein.
\newblock {The Geometry Junkyard: Dilation-Free Pla\-nar Graphs}.
\newblock Web page, 1997.
\newblock {\it http://www.ics.uci.edu/\~{}eppstein/junkyard/dilation-free/}.

\bibitem{Hillar}
C.~Hillar and D.~Rhea.
\newblock A result about the density of iterated line intersections in the
  plane.
\newblock {\em Comput. Geom. Theory Appl.}, 33(3):106--114, 2006.

\bibitem{Ismailescu}
D.~Ismailescu and R.~Radoi\v{c}i\'{c}.
\newblock A dense planar point set from iterated line intersections.
\newblock {\em Comput. Geom. Theory Appl.}, 27(3):257--267, 2004.

\bibitem{Kaz}
D.~A. Ka\v{z}dan.
\newblock Uniform distribution in the plane.
\newblock {\em Transactions of the Moscow Mathematical Society}, 14:325--332,
  1965.

\bibitem{Kutz}
R.~Klein and M.~Kutz.
\newblock The density of iterated crossing points and a gap result for
  triangulations of finite point sets.
\newblock {\em Proc. 22nd Symposium on Computational Geometry},
  Sedona, Arizona, pages 264 -- 272, ACM Press, 2006.
\end{thebibliography}
\end{document}